\def\pown{
\setlength{\textheight}{23cm}
\setlength{\textwidth}{15cm}
\setlength{\topmargin}{-2.0cm}
\setlength{\hoffset}{-0.5cm}
\setlength{\leftmargin}{-1cm}
\setlength{\rightmargin}{2.0cm}}

\documentstyle[12pt]{article}
\pown

\newcommand{\be}{\begin{equation}}
\newcommand{\ee}{\end{equation}}
\newcommand{\ba}{\begin{array}}
\newcommand{\ea}{\end{array}}
\newcommand{\beqn}{\begin{eqnarray}}
\newcommand{\eeqn}{\end{eqnarray}}
\newcommand{\beqna}{\begin{eqnarray*}}
\newcommand{\eeqna}{\end{eqnarray*}}

\newcommand{\zero}{\setcounter{equation}{0} \par}

\def\znakr{\raise1.5pt\hbox{\symb\char66\kern-2pt\char74}}
\def\znakl{\raise1.5pt\hbox{\symb\char73\kern-2pt\char67}}

\begin{document}
\title{The bicovariant differential calculus
on the $\kappa -$Poincar\`e and $\kappa -$Weyl groups}
\author{ Karol Przanowski\thanks{Supported by
\L{}\'od\'z University grant no 487} \\
Department of Field Theory \\
University of \L \'od\'z \\
ul.Pomorska 149/153, 90-236 \L \'od\'z , Poland}
\date{}
\maketitle
\begin{abstract}
The bicovariant differential calculus on four-dimen\-sio\-nal
$\kappa -$Po\-in\-care group and corresponding Lie-algebra like structure
for any metric tensor are described.
The bicovariant differential calculus on four-dimensional
$\kappa -$Weyl group and corresponding Lie-algebra like structure
for any matric tensor in the reference frame in which $g_{00}=0$
are considered.
\end{abstract}

\setcounter{section}{0}

\section{Introduction}
\zero

\def\1{$\kappa$-Poincar\'e group}
\def\2{bicovariant}
\def\3{calculus}
\def\4{differential}
\def\5{dimension}
\def\6{quantum}
\def\7{coefficient}
\def\8{operator}
\def\9{properties}
\def\0{algebra}
\def\kapo{$\kappa$-Poincar\'e}
\def\mink{$\kappa$-Min\-ko\-wski}
\def\linv{left-invariant}
\def\rinv{right-invariant}
\def\rep{representation}
\def\bic{bicrossproduct}
\def\three{three-dimensional}

Recently, considerable interest has been paid to the deformations of groups
and \0s of space-time symmetries~\cite{w1}.  In particular, an interesting
deformation of the Poincar\'e \0~\cite{w2} as well as group~\cite{w3} has
 been
introduced which depend on dimensionful deformation parameter $\kappa$; the
relevant objects are called \kapo{} \0 and \1, respectively. Their structure
was studied in some detail and many of their \9 are now well understood. The
\kapo{} \0 and group for the space-time of any \5 has been defined~\cite{w4}
, the
realizations of the \0 in terms of \4 operators acting on commutative
Min\-kowski as well as momentum spaces were given~\cite{w5};
the unitary \rep{s}
of the deformed group were found~\cite{w6}; the deformed universal covering
$ISL(2,{\cal C})$ was constructed~\cite{w7}; the bicrossproduct~\cite{w8}
 structure, both of
the \0 and group was revealed~\cite{w9}. The proof of formal duality
between \1
and \kapo{} \0 was also given, both in two~\cite{w10} as well as
in four \5s~\cite{w11}.
 One of the important problems is the construction of the \2 \4 \3 on \1.
Using an elegant approach due to Woronowicz~\cite{rach}, the \4 calculi on
four-\5al Poincar\'e group for diagonal metric tensor~\cite{diag} and
three-dimensional~\cite{tree},
 as well as on the Min\-kowski space~\cite{diag} and~\cite{min}
were constructed.

In the paper~\cite{defwayla} the $\kappa -$deformation of the
Poincar\`e algebra and group for arbitrary metric tensor has been described
and under the assumption that $g_{00}=0$ the $\kappa -$deformation of the
Weyl group as well as algebra has been constructed.

\ \

In this paper in section \ref{r1} we briefly sketch the construction
of differential calculus on the $\kappa -$Poincar\`e group for any
metric tensor $g_{\mu \nu },\ \ \mu ,\,\nu =0,...,3$.
We obtain the corresponding Lie algebra structure and prove
its equivalence to the $\kappa -$Poincar\`e algebra.
In section \ref{r2} we present the construction
of differential calculus on the $\kappa -$Weyl group for any
metric tensor $g_{\mu \nu },\ \ \mu ,\,\nu =0,...,3$ with $g_{00}=0$.
We find the corresponding Lie algebra structure and prove
its equivalence to the $\kappa -$Weyl algebra.

\section{The $\kappa -$Poincar\`e group and algebra \label{r1}}
\zero

We assume that the metric tensor $g_{\mu \nu }, (\mu \nu =0,1,...,3)$ is
represented by an arbitrary nondegenerate symetric $4\times 4$ matrix
(not necessery diagonal)
with $det(g_{\mu \nu })=1$ (in more general case in some equations we use the
 parameter $det(g)$).

The Poincar\`e group ${\cal P}$ consists of the pairs $(x,\Lambda )$,
where $x$ is a $4-$vector, $\Lambda $ is the matrix of the Lorentz
group in $4-$dimensions, with the composition law:
$$
(x^\mu ,\Lambda ^\mu _{\ \nu })*(x^{\prime \nu },
\Lambda ^{\prime \nu }_{\ \alpha }) =
(\Lambda ^\mu _{\ \nu } x^{\prime \nu }+x^\mu ,
\Lambda ^\mu _{\ \nu }\Lambda ^{\prime \nu }_{\ \alpha }) .
$$

The $\kappa -$Poincar\`e group is a Hopf $*-$algebra defined
as follows~\cite{defwayla}.
Consider the universal $*-$algebra with unity, generated by selfadjoint
elements $\Lambda ^\mu _{\ \nu },\ x^\mu $ subject to the
following relations:

\beqn
\lbrack  \Lambda ^\alpha _{\ \beta },x^\varrho  \rbrack &=& - {i\over \kappa }
((\Lambda ^\alpha _{\ 0}-\delta ^\alpha _{\ 0})
\Lambda ^\varrho _{\ \beta }+(\Lambda _{0\beta }-
 g_{0\beta })g^{\alpha \varrho }), \nonumber \\
\lbrack  x^\varrho ,x^\sigma  \rbrack  &=& {i\over \kappa }
(\delta ^\varrho _{\ 0}x^\sigma  -
\delta ^\sigma _{\ 0}x^\varrho ), \nonumber \\
\lbrack  \Lambda ^\alpha _{\ \beta }, \Lambda ^\mu _{\ \nu } \rbrack
&=& 0.  \nonumber
\eeqn

The comultiplication, antipode and counit are defined as follows:
\beqn
\Delta \Lambda ^\mu _{\ \nu } &=& \Lambda ^\mu _{\ \alpha } \otimes
\Lambda ^\alpha _{\ \nu },\nonumber \\
\Delta x^\mu  &=& \Lambda ^\mu _{\ \nu }\otimes x^\nu  +
x^\mu \otimes I, \nonumber \\
S(\Lambda ^\mu _{\ \nu }) &=& \Lambda ^{\ \mu }_\nu ,\nonumber \\
S(x^\mu )&=& -  \Lambda ^{\ \mu }_\nu  x^\nu ,\nonumber \\
\varepsilon (\Lambda ^\mu _{\ \nu })&=& \delta ^\mu _{\ \nu },\nonumber \\
\varepsilon (x^\mu )&=& 0  . \nonumber
\eeqn

In our construction of the bicovariant $*-$calculi we use the Woronowicz
theory,~\cite{rach}.
First we construct the right ad-invariant ideal ${\cal R}$ in $ker\  \varepsilon,
\   (ad({\cal R})={\cal R}\otimes {\cal P}_\kappa )$.
The adjoint action of the group is defined as follows:
$$
ad(a)=\sum _{k} b_k\otimes S(a_k)c_k ,
$$
where
$$
(\Delta \otimes I)\Delta (a) = (I\otimes \Delta )\Delta (a) = \sum _{k}
a_k\otimes b_k\otimes c_k .
$$

In order to obtain the ideal ${\cal R}$ we put:
\beqn
\Delta ^\mu _{\ \nu} &=& \Lambda ^\mu _{\ \nu}
- \delta ^\mu _{\ \nu},   \nonumber \\
\Delta ^{\mu \ \alpha }_{\ \nu} &=& x^\alpha \Delta ^\mu _{\ \nu}
+ {i\over \kappa }(g_{\nu 0}\Delta ^{\alpha \mu }
- \delta ^\mu _{\ 0}\Delta ^{\alpha }_{\ \nu }),   \nonumber \\
x^{\alpha \beta } &=& x^\alpha  x^\beta  + {i\over \kappa }
(g^{\alpha \beta } x_0 - \delta ^\alpha _{\ 0} x^\beta ), \nonumber \\
{\tilde \Delta} ^{\mu \ \alpha }_{\ \nu} &=& \Delta ^{\mu \ \alpha }_{\ \nu}
- \fbox{$ {1\over 6} \varepsilon ^{\mu \ \alpha \gamma }_{\ \nu }
\varepsilon _{\gamma \varrho \sigma \delta }
\Delta ^{\varrho \sigma \delta }$}, \nonumber \\
{\tilde x}^{\alpha \beta } &=& x^{\alpha \beta }
- \fbox{$ {1\over 4} g^{\alpha \beta } x^\mu _{\ \mu } $} , \nonumber \\
T &=& lin\{\ {\tilde \Delta }^{\mu \ \alpha }_{\ \nu},
\ \Delta ^\alpha _{\ \beta }\Delta ^\mu _{\ \nu },
\ {\tilde x}^{\alpha \beta }\   \}.  \nonumber
\eeqn

One can prove that the linear set $T$ is ad-invariant.
Let ${\cal R}$ be a right ideal of ${\cal P_\kappa }$
generated by elements of $T$.
Then the following theorem holds:

{\bf Theorem 1.} ${\cal R}$ has the following properties:

(i)\ ${\cal R}$ is ad-invariant,

(ii)\ for any $a\in {\cal R}, \ S(a)^* \in {\cal R}$,

(iii)\ $ker\ \varepsilon /{\cal R}$ is spanned by the following
elements:

$$
x^\mu ;\ \ \ \Delta ^\mu _{\ \nu },\ \mu <\nu ;\ \ \
\varphi = x^\mu _{\ \mu };
$$

$$
\varphi _\mu = \varepsilon _{\mu \varrho \sigma \delta }
\Delta ^{\varrho \sigma \delta }.
$$

Note only that:

(a)\ $\Delta ^{\mu \ \alpha }_{\ \nu },\ x^{\mu \nu }$ are "improved"
generators of $(ker\  \varepsilon )^2$ which form the (not completely
reducible) multiplet under adjoint action of ${\cal P}_\kappa $,

(b)\ if $g_{00}\ne 0$ then the ideal generated by $\Delta ^\alpha _{\ \beta }
 \Delta ^\mu _{\ \nu },\  \Delta ^{\mu \ \alpha }_{\ \nu }$,
 $x^{\mu \nu }$ equals $ker\ \varepsilon $. In order to obtain reasonable
 (in the sesne that it contains all differentials $dx^\mu ,\
 d\Lambda ^\mu _{\ \nu }$) calculus we have subtracted the trace
 of $x^{\mu \nu }$ and completely antisymmetric part of
 $\Delta ^{\mu \ \alpha }_{\ \nu }$.

(c)\ it is easy to conclude from (iii) that our calculus is fifteen-dimensional.

Now, following Woronowicz's paper~\cite{rach} we find
the basis of the space of the left-invariant 1-forms:

\beqn
\omega ^\mu  &=& \pi r^{-1}(I\otimes x^\mu ) =
\Lambda ^{\ \mu }_\nu dx^\nu , \nonumber \\
\omega ^\mu _{\ \nu } &=& \pi r^{-1}(I\otimes \Delta ^\mu _{\ \nu })
= \Lambda ^{\ \mu }_\alpha d\Lambda ^\alpha _{\ \nu }, \nonumber \\
\Omega  &=& \pi r^{-1}(I\otimes \varphi ) = d\varphi
-2 x_\alpha dx^\alpha , \nonumber \\
\Omega _\mu  &=& \pi r^{-1}(I\otimes \varphi _\mu ) =
\varepsilon _{\mu \nu \alpha \beta }\Lambda ^{\ \nu }_\delta
\omega ^\beta \Lambda ^{\delta \alpha }
-{2i\over \kappa }\varepsilon _{0\mu \nu \beta }
\omega ^{\nu \beta }. \nonumber
\eeqn

The commutations rules between the invariant forms and generators of
${\cal P}_\kappa $ read:

\beqn
\lbrack \Lambda ^\mu _{\ \nu },\omega ^\alpha   \rbrack  &=&
-{i\over \kappa }(g_{\nu 0}\Lambda ^\mu _{\ \tau }\omega ^{\tau \alpha }
+\Lambda ^\mu _{\ 0}\omega ^\alpha _{\ \nu })
-{1\over 6}\varepsilon ^{\tau \ \alpha \gamma }_{\ \nu }
\Lambda ^\mu _{\ \tau }\Omega _\gamma , \nonumber \\
\lbrack x^\mu ,\omega ^\alpha   \rbrack  &=&
-{1\over 4} \Lambda ^{\mu \alpha }\Omega
+{i\over \kappa } (\Lambda ^{\mu \alpha }\omega _0
-\delta ^\alpha _{\ 0}\Lambda ^\mu _{\ \nu }\omega ^\nu ), \nonumber \\
\lbrack \Lambda ^\alpha _{\ \beta },\omega ^\mu _{\ \nu }  \rbrack
 &=& 0, \nonumber \\
\lbrack  x^\alpha ,\omega ^\mu _{\ \nu } \rbrack  &=&
-{i\over \kappa }(\delta ^\mu _{\ 0}\Lambda ^\alpha _{\ \beta }
\omega ^\beta _{\ \nu }
+g_{\nu 0}\Lambda ^\alpha _{\ \beta }\omega ^{\mu \beta }
-\Lambda ^\alpha _{\ \nu }\omega ^\mu _{\ 0}
-\Lambda ^{\alpha \mu }\omega _{0\nu })
-{1\over 6}\varepsilon ^{\mu \ \beta \gamma }_{\ \nu }
 \Lambda ^\alpha _{\ \beta }\Omega _\gamma , \nonumber \\
\lbrack \Lambda ^\mu _{\ \nu },\Omega   \rbrack  &=&
{4\over \kappa ^2}g_{00}\Lambda ^{\mu \tau }
\omega _{\tau \nu }, \nonumber \\
\lbrack x^\alpha  ,\Omega  \rbrack  &=&
{4\over \kappa ^2}g_{00}\Lambda ^\alpha _{\ \beta }
\omega ^\beta , \nonumber \\
\lbrack  \Lambda ^\alpha _{\ \beta },\Omega _\mu   \rbrack  &=& 0, \nonumber \\
\lbrack x^\alpha ,\Omega _\mu   \rbrack  &=&
{3\over \kappa ^2}g_{00}\varepsilon _{\mu \nu \tau \beta }
\Lambda ^{\alpha \nu }\omega ^{\tau \beta }
+{i\over \kappa }(\Lambda ^\alpha _{\ \mu }\Omega _0
-g_{0\mu }\Lambda ^{\alpha \beta }\Omega _\beta ). \nonumber
\eeqn

The corresponding basis of the space of the right-invariant 1-forms is:

\beqn
\eta ^\mu  &=&
-\omega ^\alpha _{\ \beta }\Lambda ^{\ \beta }_\nu
x^\nu \Lambda ^\mu _{\ \alpha }
+\omega ^\alpha  \Lambda ^\mu _{\ \alpha },\nonumber  \\
\eta ^\mu _{\ \nu } &=& \omega ^\alpha _{\ \beta }
\Lambda ^\mu _{\ \alpha }\Lambda ^{\ \beta }_\nu  ,\nonumber  \\
\theta  &=& \Omega  ,\nonumber  \\
\theta _\mu  &=& \Omega _\nu  \Lambda ^{\ \nu }_\mu  .\nonumber
\eeqn

%
%
%

This concludes the description of bimodule $\Gamma $ of 1-forms on
${\cal P}_\kappa $. The external algebra can be now constructed
as follows~\cite{rach}. On $\Gamma ^{\otimes 2}$ we define a bimodule
homomorpfism $\sigma $ such that
$$
\sigma (\omega \otimes _{{\cal P}_\kappa } \eta) =
\eta \otimes _{{\cal P}_\kappa } \omega ,
$$
for any left-invariant $\omega \in \Gamma $ and any right-invariant
$\eta \in \Gamma $. Then by the definition

$$
\Gamma ^{\wedge 2} = {\Gamma ^{\otimes 2}\over ker(I-\sigma )}.
$$

Finally, after a long analysis we obtain the following set of relations:

\beqn
 &\ & \Omega \wedge \Omega  = 0, \nonumber \\
 &\ & \omega ^\mu _{\ \nu }\wedge \omega ^\alpha _{\ \beta }
 +\omega ^\alpha _{\ \beta }\wedge \omega ^\mu _{\ \nu } = 0, \nonumber \\
 &\ & \Omega _\alpha \wedge \Omega _\beta
 +\Omega _\beta \wedge \Omega _\alpha  = 0, \nonumber \\
 &\ & \Omega _\alpha \wedge \omega ^\mu _{\ \nu }
 +\omega ^\mu _{\ \nu }\wedge \Omega _\alpha  = 0, \nonumber \\
 &\ & \Omega \wedge \omega ^\mu _{\ \nu }
 +\omega ^\mu _{\ \nu }\wedge \Omega
 -{4\over \kappa ^2}g_{00}\omega _{\beta \nu }
 \wedge \omega ^{\beta \mu } = 0, \nonumber \\
 &\ & \Omega _\mu  \wedge \Omega
 +\Omega \wedge \Omega _\mu  -{4\over \kappa ^2}g_{00}
 \Omega _\tau \wedge \omega ^\tau _{\ \mu }
 -{4\over \kappa ^2}g_{00}{\cal X}_\mu  = 0, \nonumber \\
 &\ & \omega ^\mu \wedge \Omega
 +\Omega \wedge \omega ^\mu
 -{4\over \kappa ^2}g_{00}\omega ^\mu _{\ \beta }
 \wedge \omega ^\beta = 0,\nonumber \\
 &\ & \omega ^\alpha \wedge \omega ^\mu
 +\omega ^\mu \wedge \omega ^\alpha
 +{i\over \kappa }(\delta ^\alpha _{\ 0}\omega ^\mu _{\ \beta }
 \wedge \omega ^\beta
 +\delta ^\mu _{\ 0} \omega ^\alpha _{\ \beta }
 \wedge \omega ^\beta ) = 0, \nonumber \\
 &\ & \omega ^\alpha \wedge \omega ^{\mu \nu }
 +\omega ^{\mu \nu }\wedge \omega ^\alpha
 +{i\over \kappa }(\delta ^\mu _{\ 0}\omega ^{\tau \nu }
 \wedge \omega ^{\ \alpha }_\tau
 +\delta ^\nu _{\ 0}\omega ^\mu _{\ \tau }\wedge
 \omega ^{\tau \alpha }
 +\omega ^{\ \nu }_0 \wedge\omega ^{\alpha \mu }
 +\omega ^\mu _{\ 0}\wedge \omega ^{\alpha \nu })  \nonumber \\
&\ & -{1\over 6}\varepsilon ^{\alpha \mu \nu \sigma }
 {\cal X}_\sigma = 0,  \nonumber \\
 &\ & \omega ^\mu \wedge \Omega ^\alpha
 +\Omega ^\alpha \wedge \omega ^\mu
 +{i\over \kappa }(\delta ^\alpha _{\ 0}\Omega _\tau \wedge
 \omega ^{\tau \mu }
 +\Omega _0 \wedge \omega ^{\mu \alpha }   \nonumber \\
&\ & +{1\over 12}\varepsilon ^{\alpha \mu \tau \sigma }
 \Omega _\tau \wedge \Omega _\sigma
 -{3\over 2\kappa ^2}g_{00}\varepsilon ^{\alpha \sigma \tau \beta }
 \omega ^\mu _{\ \sigma }\wedge \omega _{\tau \beta }
 -{1\over 4}g^{\alpha \mu }{\cal Y} = 0, \nonumber
\eeqn

where we introduce the following 2-forms:
\beqn
&\ & {\cal Y} = \omega _\alpha \wedge \Omega ^\alpha
+\Omega ^\alpha \wedge \omega _\alpha
+{i\over \kappa }\Omega _\tau \wedge \omega ^\tau _{\ 0}, \nonumber \\
&\ & {\cal X}_\sigma = \varepsilon _{\sigma \alpha \mu \nu }
(\omega ^\alpha \wedge \omega ^{\mu \nu }
+\omega ^{\mu \nu }\wedge \omega ^\alpha
+{2i\over \kappa }(\delta ^\mu _{\ 0}
\omega ^{\tau \nu }\wedge \omega ^{\ \alpha }_\tau
+\omega ^{\ \nu }_0 \wedge \omega ^{\alpha \mu })). \nonumber
\eeqn

The basis of $\Gamma ^{\wedge 2}$ consists of the following elements:

\beqn
 &\ & \omega ^{\alpha \beta }\wedge \omega ^{\mu \nu };
 \ \ \ \ {\rm for}\ \alpha <\beta ,\ \mu < \nu,\ (\alpha \beta )
 \ne (\mu \nu ),\ \alpha < \mu ; \nonumber \\
 &\ & \Omega ^\mu \wedge \Omega ^\nu ,
 \ \omega ^\mu \wedge \Omega ^\nu ,
 \ \omega ^{\mu \nu }\wedge \Omega ^\alpha ,
 \ \omega ^{\mu \nu }\wedge \Omega ,
 \ \omega ^\mu \wedge \omega ^\nu ,
 \ \omega ^{\mu \nu }\wedge \omega ^\alpha ;
 \ \ \ \ {\rm for}\ \mu <\nu ; \nonumber \\
 &\ &  \omega ^\mu \wedge \Omega ,
 \ \Omega ^\mu \wedge \Omega ,
 \ {\cal Y},
 \ {\cal X}^\mu . \nonumber
\eeqn

Thus, there are five more elements than it is generically expected.
The Cartan-Maurer equations have the following form:
\beqn
d\omega ^\mu _{\ \nu } &=& \omega ^{\ \mu }_\tau
\wedge \omega ^\tau _{\ \nu }, \nonumber \\
d\omega ^\mu  &=&  \omega ^{\ \mu }_\tau \wedge \omega ^\tau , \nonumber \\
d\Omega  &=& 0, \nonumber \\
d\Omega _\mu  &=& -\omega ^{\ \alpha }_\mu \wedge
\Omega _\alpha  - {\cal X}_\mu . \nonumber
\eeqn

To obtain the quantum Lie algebra we introduce the left-invariant fields,
defined by the formula:
\beqn
da = {1\over 2}(\chi _{\mu \nu }*a)\omega ^{\mu \nu }
+(\chi _\mu *a)\omega ^\mu
+(\chi *a)\Omega
+(\lambda _\mu *a)\Omega ^\mu , \label{p_pole}
\eeqn
where, for any linear functional $\varphi $ on ${\cal P}_\kappa $,
$$
\varphi *a = (I\otimes \varphi )\Delta (a).
$$
The product of two functional $\varphi _1,\ \varphi _2$ is defined
by the duality relation:
$$
\varphi _1 * \varphi _2 (a) = (\varphi _1 \otimes \varphi _2)
\Delta (a). 
$$
Finally, we apply the external derivative to both sides of
eq.(\ref{p_pole}). Using the fact that $d^2 a=0$ and equating to zero
the coefficients of the basis elements of $\Gamma ^{\wedge 2}$
we find the quantum Lie algebra:
\beqn
\lambda^\mu  \chi _\mu &=& 0, \nonumber \\
\lbrack  \chi_{\mu \nu }, \chi  \rbrack  &=& 0, \nonumber \\
\lambda _\mu ({4\over \kappa ^2}g_{00}\chi -1)&=& {1\over 12}
\varepsilon ^{\ \alpha \beta \nu }_\mu \chi _\alpha
\chi _{\beta \nu }, \nonumber \\
\lbrack \chi _{\mu \nu },\chi _\alpha  \rbrack &=&
(1+{i\over \kappa }\chi _0 -{4\over \kappa ^2}g_{00}\chi )
(\chi _\mu g_{\alpha \nu }-\chi _\nu g_{\alpha \mu }) \nonumber \\
\lbrack \lambda _\mu ,\chi   \rbrack &=& 0, \nonumber \\
\lbrack \chi _\alpha ,\chi _\mu   \rbrack &=& 0  \nonumber \\
\lbrack \chi ,\chi _\mu   \rbrack &=& 0, \nonumber \\
\lbrack \chi _{\mu \nu },\lambda _\alpha   \rbrack &=&
(\lambda _\mu g_{\alpha \nu }-\lambda _\nu g_{\alpha \mu })
(1-{4\over \kappa ^2}g_{00}\chi )    \nonumber \\
&\ & -{i\over \kappa }\lambda _0 (\chi _\nu g_{\mu \alpha }
-\chi _\mu g_{\nu \alpha })
-{i\over \kappa }g_{0\alpha }(\lambda _\nu \chi _\mu
-\lambda _\mu \chi _\nu ),  \nonumber  \\
\lbrack \lambda _\alpha ,\chi _\mu   \rbrack  &=& 0, \nonumber \\
\lbrack \lambda _\alpha ,\lambda _\mu   \rbrack &=& -{i\over 6}
\varepsilon ^{\ \ \ \sigma \delta }_{\mu \alpha }\lambda _\sigma
\chi _\delta,  \nonumber \\
\lbrack \chi _{\alpha \beta },\chi _{\mu \nu } \rbrack &=&
(1-{4\over \kappa ^2}g_{00}\chi )
(\chi _{\beta \mu }g_{\nu \alpha }
+\chi _{\alpha \nu }g_{\mu \beta }
-\chi _{\beta \nu }g_{\mu \alpha }
-\chi _{\alpha \mu }g_{\nu \beta }) \nonumber \\
&\ & +{i\over \kappa }(\chi _\alpha (
\chi _{\mu 0}g_{\nu \beta }-\chi _{\nu 0}g_{\mu \beta })
+\chi _\mu (\chi _{\beta 0}g_{\alpha \nu }
-\chi _{\alpha 0}g_{\nu \beta }) \nonumber \\
&\ & +\chi _\nu (\chi _{\alpha 0}g_{\mu \beta }
-\chi _{\beta 0}g_{\alpha \mu })
+\chi _\beta (\chi _{\nu 0}g_{\alpha \mu }
-\chi _{\mu 0}g_{\alpha \nu })) \nonumber \\
&\ & +{i\over \kappa }(\chi _\beta (\chi _{\mu \alpha }g_{0\nu }
-\chi _{\nu \alpha }g_{0\mu })
+\chi _\nu (\chi _{\beta \mu }g_{0\alpha }
-\chi _{\alpha \mu }g_{0\beta }) \nonumber \\
&\ & +\chi _\mu (\chi _{\alpha \nu }g_{0\beta }
-\chi _{\beta \nu }g_{0\alpha })
+\chi _\alpha (\chi _{\nu \beta }g_{0\mu }
-\chi _{\mu \beta }g_{0\nu }))  \nonumber \\
&\ & +{3\over \kappa ^2}g_{00}\lambda _\sigma
(\chi _\beta  \varepsilon ^\sigma _{\ \alpha \mu \nu }
-\chi _\alpha  \varepsilon ^\sigma _{\ \beta \mu \nu }
+\chi _\mu  \varepsilon ^\sigma _{\ \nu \alpha \beta }
-\chi _\nu  \varepsilon ^\sigma _{\ \mu \alpha \beta }
). \nonumber
\eeqn

Having our quantum Lie algebra constructed, we can now pose the question
what the relation is between our functionals and the elements of the
$\kappa -$Poincare algebra $\tilde {\cal P}_\kappa $.

The $\kappa -$Poincare
algebra is defined as follows~\cite{defwayla},~\cite{w11}:

\noindent  The commutation rules:
\beqn
\lbrack  M^{ij},P_0 \rbrack  &=& 0 ,\nonumber \\
\lbrack M^{ij},P_k  \rbrack  &=& i\kappa (\delta ^j_{\ k}
g^{0i}-\delta ^i_{\ k}g^{0j})(1-e^{-{P_0\over \kappa }})
+ i(\delta ^j_{\ k}g^{is}-\delta ^i_{\ k}g^{js})P_s ,\nonumber \\
\lbrack M^{i0},P_0  \rbrack  &=& i\kappa g^{i0}
(1-e^{-{p_0\over \kappa }})+ig^{ik}P_k ,\nonumber \\
\lbrack M^{i0},P_k  \rbrack  &=& -i{\kappa \over 2}g^{00}
\delta ^i_{\ k}(1-e^{-2{P_0\over \kappa }})-
i\delta ^i_{\ k}g^{0s}P_s e^{-{P_0\over \kappa }}+ ,\nonumber \\
 &\ & +ig^{0i}P_k (e^{-{P_0\over \kappa }}-1)+
{i\over 2\kappa }\delta ^i_{\ k}g^{rs}P_r P_s -
{i\over \kappa }g^{is} P_s P_k ,\nonumber \\
\lbrack P_\mu ,P_\nu   \rbrack  &=& 0 ,\nonumber \\
\lbrack M^{\mu \nu },M^{\lambda \sigma }  \rbrack
&=& i(g^{\mu \sigma }M^{\nu \lambda }-g^{\nu \sigma }M^{\mu \lambda }
+g^{\nu \lambda }M^{\mu \sigma }-g^{\mu \lambda }M^{\nu \sigma }) .\nonumber
\eeqn
The coproducts, counit and antipode:
\beqn
\Delta P_0 &=& I\otimes P_0 + P_0 \otimes I ,\nonumber \\
\Delta P_k &=& P_k \otimes e^{-{P_0\over \kappa }}+
I\otimes P_k ,\nonumber \\
\Delta M^{ij} &=& M^{ij}\otimes I + I\otimes M^{ij} ,\nonumber \\
\Delta M^{i0} &=& I\otimes M^{i0} + M^{i0}\otimes
e^{-{P_0\over \kappa }}- {1\over \kappa }M^{ij}\otimes P_j ,\nonumber \\
\varepsilon (M^{\mu \nu }) &=& 0; \ \ \
\varepsilon (P_\nu ) = 0;\ \ \
\varepsilon (D) =0, \nonumber \\
S(P_0) &=& -P_0 ,\nonumber \\
S(P_i) &=& -e^{P_0\over \kappa }P_i ,\nonumber \\
S(M^{ij}) &=& -M^{ij} ,\nonumber \\
S(M^{i0}) &=& -e^{P_0\over \kappa }(M^{i0}+
{1\over \kappa }M^{ij}P_j) ,\nonumber
\eeqn
where $i,j,k = 1,2,3$.

Using, on the one hand, the properties of the left-invariant fields
described by Woronowicz~\cite{rach} and on the other hand,
the duality relations ${\cal P}_\kappa \iff \tilde {\cal P}_\kappa $
established in~\cite{w11}, one can prove that the following substitutions:

\beqn
\chi _0  &=& -i(\kappa (e^{P_0\over \kappa }-1)
-{g_{00}\over 2\kappa }M^2), \nonumber \\
\chi _i  &=& -i(e^{P_0\over \kappa }P_i
-g_{i0}{M^2 \over 2\kappa }), \nonumber \\
\chi   &=& -{M^2\over 8}, \nonumber \\
\chi _{ij} &=& -i (1+{i\over \kappa }\chi _0
-{4\over \kappa ^2}g_{00}\chi )M_{ij}
-{1\over \kappa }(\chi _0 M_{ij} +\chi _i M_{j0}- \chi _j M_{i0}), \nonumber \\
\chi _{i0} &=& -i (1+{i\over \kappa }\chi _0
-{4\over \kappa ^2}g_{00}\chi )M_{i0}, \nonumber
\eeqn
where
$$
M^2 = g^{00}(2\kappa \sinh({P_0\over 2\kappa }))^2
+4\kappa g^{0l} P_l e^{P_0\over 2\kappa } \sinh({P_0\over 2\kappa })
+g^{rs}P_r e^{P_0\over 2\kappa }P_s e^{P_0\over 2\kappa },
$$

reproduce the algebra and coalgebra structure of our quantum Lie algebra.

Now, it is easy to see that the vectorfields $\chi $ and $\lambda _\mu $
are proportional to first casimir operator and deformed Paul-Luba\`nski
invariant.

\section{The $\kappa -$Poincar\`e group and algebra in case $g_{00}=0$}
\zero

After constructing this calculus, prof.J.Lukierski
suggested me that the assumption that $g_{00}=0$ should simplify
this calculus. This problem in case of
the differential calculus on $\kappa -$Min\-kowski space,
is discussed in~\cite{luk}.It appears that in the case
of differential calculus on
the $\kappa -$Poincar\`e group, under the assumption $g_{00}=0$ we can
obtain the differential calculus whose dimension
is equal to the dimension of the classical differential calculus.
In this case the ideal generated by the elements:
$\Delta ^\alpha _{\ \beta }
\Delta ^\mu _{\ \nu },\ \Delta ^{\mu \ \alpha }_{\ \nu },\
x^{\mu \nu }$
is adjoint invariant and is not equal to $ker\ \varepsilon $.
This ideal gives us the ten-dimensional calculus and the dimension of
square exterior power of our differential calculus is $10\choose 2$ dimensional.
But if $g_{00}\ne 0$ the situation dramatically changes, and we obtain
the sexteen-dimensional differential calculus and the exterior power
of our calculus need additional five differential forms.

To obtain the ideal ${\cal R}$ we put:
\beqn
T &=& lin\{\ \Delta ^{\mu \ \alpha }_{\ \nu},
\ \Delta ^\alpha _{\ \beta }\Delta ^\mu _{\ \nu },
\ x^{\alpha \beta }\   \}.  \nonumber
\eeqn
In the definition of $T$ we not subtract the trace of $x^{\mu \nu }$
and completely antisymmetric part of
 $\Delta ^{\mu \ \alpha }_{\ \nu }$.

One can prove that the linear set $T$ is ad-invariant.
Let ${\cal R}$ be a right ideal of ${\cal P_\kappa }$
generated by elements of $T$.
Then the following theorem holds:

{\bf Theorem 2.} ${\cal R}$ has the following properties:

(i)\ ${\cal R}$ is ad-invariant,

(ii)\ for any $a\in {\cal R}, \ S(a)^* \in {\cal R}$,

(iii)\ $ker\ \varepsilon /{\cal R}$ is spanned by the following
elements:

$$
x^\mu ;\ \ \ \Delta ^\mu _{\ \nu },\ \mu <\nu ;\ \ \
$$

The basis of the space of the left-invariant 1-forms read:

\beqn
\omega ^\mu  &=& \pi r^{-1}(I\otimes x^\mu ) =
\Lambda ^{\ \mu }_\nu dx^\nu , \nonumber \\
\omega ^\mu _{\ \nu } &=& \pi r^{-1}(I\otimes \Delta ^\mu _{\ \nu })
= \Lambda ^{\ \mu }_\alpha d\Lambda ^\alpha _{\ \nu }. \label{fp}
\eeqn

The commutations rules between the invariant forms and generators of
${\cal P}_\kappa $ read:

\beqn
\lbrack \Lambda ^\mu _{\ \nu },\omega ^\alpha   \rbrack  &=&
-{i\over \kappa }(g_{\nu 0}\Lambda ^\mu _{\ \tau }\omega ^{\tau \alpha }
+\Lambda ^\mu _{\ 0}\omega ^\alpha _{\ \nu }) , \nonumber \\
\lbrack x^\mu ,\omega ^\alpha   \rbrack  &=&
{i\over \kappa } (\Lambda ^{\mu \alpha }\omega _0
-\delta ^\alpha _{\ 0}\Lambda ^\mu _{\ \nu }\omega ^\nu ), \nonumber \\
\lbrack \Lambda ^\alpha _{\ \beta },\omega ^\mu _{\ \nu }  \rbrack
 &=& 0, \nonumber \\
\lbrack  x^\alpha ,\omega ^\mu _{\ \nu } \rbrack  &=&
-{i\over \kappa }(\delta ^\mu _{\ 0}\Lambda ^\alpha _{\ \beta }
\omega ^\beta _{\ \nu }
+g_{\nu 0}\Lambda ^\alpha _{\ \beta }\omega ^{\mu \beta } \nonumber \\
&\ &-\Lambda ^\alpha _{\ \nu }\omega ^\mu _{\ 0}
-\Lambda ^{\alpha \mu }\omega _{0\nu }) . \label{kp}
\eeqn

The external algebra read:

\beqn
 &\ & \omega ^\mu _{\ \nu }\wedge \omega ^\alpha _{\ \beta }
 +\omega ^\alpha _{\ \beta }\wedge \omega ^\mu _{\ \nu } = 0, \nonumber \\
 &\ & \omega ^\alpha \wedge \omega ^\mu
 +\omega ^\mu \wedge \omega ^\alpha
 +{i\over \kappa }(\delta ^\alpha _{\ 0}\omega ^\mu _{\ \beta }
 \wedge \omega ^\beta
 +\delta ^\mu _{\ 0} \omega ^\alpha _{\ \beta }
 \wedge \omega ^\beta ) = 0, \nonumber \\
 &\ & \omega ^\alpha \wedge \omega ^{\mu \nu }
 +\omega ^{\mu \nu }\wedge \omega ^\alpha
 +{i\over \kappa }(\delta ^\mu _{\ 0}\omega ^{\tau \nu }
 \wedge \omega ^{\ \alpha }_\tau
 +\delta ^\nu _{\ 0}\omega ^\mu _{\ \tau }\wedge
 \omega ^{\tau \alpha }   \nonumber \\
&\ & +\omega ^{\ \nu }_0 \wedge\omega ^{\alpha \mu }
 +\omega ^\mu _{\ 0}\wedge \omega ^{\alpha \nu }) =0.  \label{zp}
\eeqn

The basis of $\Gamma ^{\wedge 2}$ consists of the following elements:

\beqn
 &\ & \omega ^{\alpha \beta }\wedge \omega ^{\mu \nu };
 \ \ \ \ {\rm for}\ \alpha <\beta ,\ \mu < \nu,\ (\alpha \beta )
 \ne (\mu \nu ),\ \alpha < \mu ; \nonumber \\
 &\ & \omega ^\mu \wedge \omega ^\nu ,
 \ \omega ^{\mu \nu }\wedge \omega ^\alpha ;
 \ \ \ \ {\rm for}\ \mu <\nu . \nonumber
\eeqn

The Cartan-Maurer equations have the following form:
\beqn
d\omega ^\mu _{\ \nu } &=& \omega ^{\ \mu }_\tau
\wedge \omega ^\tau _{\ \nu }, \nonumber \\
d\omega ^\mu  &=&  \omega ^{\ \mu }_\tau \wedge \omega ^\tau . \label{mp}
\eeqn

Note that this calculus can be obtained from the previous one by putting
$g_{00}=0,\ \Omega =\Omega _\mu ={\cal X}_\mu ={\cal Y}=0$.

Introduceing the left-invariant fields,
defined by the formula:
$$
da = {1\over 2}(\chi _{\mu \nu }*a)\omega ^{\mu \nu }
+(\chi _\mu *a)\omega ^\mu ,
$$

we find the quantum Lie algebra:
\beqn
\lbrack \chi _{\mu \nu },\chi _\alpha  \rbrack &=&
(1+{i\over \kappa }\chi _0 )
(\chi _\mu g_{\alpha \nu }-\chi _\nu g_{\alpha \mu }) \nonumber \\
\lbrack \chi _\alpha ,\chi _\mu   \rbrack &=& 0  \nonumber \\
\lbrack \chi _{\alpha \beta },\chi _{\mu \nu } \rbrack &=&
(\chi _{\beta \mu }g_{\nu \alpha }
+\chi _{\alpha \nu }g_{\mu \beta }
-\chi _{\beta \nu }g_{\mu \alpha }
-\chi _{\alpha \mu }g_{\nu \beta }) \nonumber \\
 +{i\over \kappa }(\chi _\alpha (
\chi _{\mu 0}g_{\nu \beta }&-&\chi _{\nu 0}g_{\mu \beta })
+\chi _\mu (\chi _{\beta 0}g_{\alpha \nu }
-\chi _{\alpha 0}g_{\nu \beta }) \nonumber \\
 +\chi _\nu (\chi _{\alpha 0}g_{\mu \beta }
&-&\chi _{\beta 0}g_{\alpha \mu })
+\chi _\beta (\chi _{\nu 0}g_{\alpha \mu }
-\chi _{\mu 0}g_{\alpha \nu })) \nonumber \\
 +{i\over \kappa }(\chi _\beta (\chi _{\mu \alpha }g_{0\nu }
&-&\chi _{\nu \alpha }g_{0\mu })
+\chi _\nu (\chi _{\beta \mu }g_{0\alpha }
-\chi _{\alpha \mu }g_{0\beta }) \nonumber \\
 +\chi _\mu (\chi _{\alpha \nu }g_{0\beta }
&-&\chi _{\beta \nu }g_{0\alpha })
+\chi _\alpha (\chi _{\nu \beta }g_{0\mu }
-\chi _{\mu \beta }g_{0\nu }) . \nonumber
\eeqn

One can prove that the following substitutions:

\beqn
\chi _0  &=& -i\kappa (e^{P_0\over \kappa }-1), \nonumber \\
\chi _i  &=& -i(e^{P_0\over \kappa }P_i
-g_{i0}{M^2 \over 2\kappa }), \nonumber \\
\chi _{ij} &=& -i (1+{i\over \kappa }\chi _0 )M_{ij}
-{1\over \kappa }(\chi _0 M_{ij} +\chi _i M_{j0}- \chi _j M_{i0}), \nonumber \\
\chi _{i0} &=& -i (1+{i\over \kappa }\chi _0 )M_{i0}, \nonumber
\eeqn
where
\beqn
M^2 &=& g^{00}(2\kappa \sinh({P_0\over 2\kappa }))^2
+4\kappa g^{0l} P_l e^{P_0\over 2\kappa }
\sinh({P_0\over 2\kappa }) \nonumber \\
&\ &+g^{rs}P_r e^{P_0\over 2\kappa }P_s e^{P_0\over 2\kappa },  \nonumber
\eeqn

reproduce the algebra and coalgebra structure of our quantum Lie algebra.

\section{The $\kappa -$Weyl group and algebra \label{r2}}
\zero

The Weyl group ${\cal W}$ consists of the triples $(x,\Lambda ,e^b)$,
where $x$ is a $4-$vector, $\Lambda $ is the matrix of the Lorentz
group in $4-$dimensions and $b\in R$, with the composition law:
$$
(x^\mu ,\Lambda ^\mu _{\ \nu },e^b)*(x^{\prime \nu },
\Lambda ^{\prime \nu }_{\ \alpha },e^{b^\prime }) =
(\Lambda ^\mu _{\ \nu }e^b x^{\prime \nu }+x^\mu ,
\Lambda ^\mu _{\ \nu }\Lambda ^{\prime \nu }_{\ \alpha },
e^b e^{b^\prime }) .
$$

In~\cite{defwayla}, under assumption that $g_{00}=0$ the
$\kappa -$deformation of the Weyl group ${\cal W}_\kappa $ was constructed.
This $\kappa -$Weyl group ${\cal W}_\kappa $ is a Hopf $*-$algebra defined
as follows~\cite{defwayla}. Consider the universal $*-$algebra with
unity generated by self adjoint elements $\Lambda ^\mu _{\ \nu },\
x^\mu $ and $e^b$ subject to the following relations:

\beqn
\lbrack  \Lambda ^\alpha _{\ \beta },x^\varrho  \rbrack &=& - {i\over \kappa }
((e^b \Lambda ^\alpha _{\ 0}-\delta ^\alpha _{\ 0})
\Lambda ^\varrho _{\ \beta }+(\Lambda _{0\beta }-
e^b g_{0\beta })g^{\alpha \varrho }), \nonumber \\
\lbrack  x^\varrho ,x^\sigma  \rbrack  &=& {i\over \kappa }
(\delta ^\varrho _{\ 0}x^\sigma  -
\delta ^\sigma _{\ 0}x^\varrho ), \nonumber \\
\lbrack  \Lambda ^\alpha _{\ \beta }, \Lambda ^\mu _{\ \nu } \rbrack
&=& 0,  \nonumber \\
\lbrack \Lambda ^\mu _{\ \nu },e^b \rbrack  &=& 0, \nonumber \\
\lbrack x^\mu ,e^b \rbrack  &=& 0 . \nonumber
\eeqn

The coproduct, antipode and counit are defined as follows:
\beqn
\Delta \Lambda ^\mu _{\ \nu } &=& \Lambda ^\mu _{\ \alpha } \otimes
\Lambda ^\alpha _{\ \nu },\nonumber \\
\Delta x^\mu  &=& e^b \Lambda ^\mu _{\ \nu }\otimes x^\nu  +
x^\mu \otimes I,   \nonumber \\
\Delta e^b &=& e^b\otimes e^b  , \nonumber \\
S(\Lambda ^\mu _{\ \nu }) &=& \Lambda ^{\ \mu }_\nu ,\nonumber \\
S(x^\mu )&=& - e^{-b} \Lambda ^{\ \mu }_\nu  x^\nu ,\nonumber \\
S(e^b) &=& e^{-b}, \nonumber \\
\varepsilon (\Lambda ^\mu _{\ \nu })&=& \delta ^\mu _{\ \nu },\nonumber \\
\varepsilon (x^\mu )&=& 0,  \nonumber \\
\varepsilon (e^b) &=& 1 . \nonumber
\eeqn

In order to obtain the bicovariant $*-$calculi we go along the same
lines as in the $\kappa -$Poincare group.

Let ${\cal R}$ be a right ideal of ${\cal W_\kappa }$
generated by the following elements:
$(e^b - I)^2,\ \Delta ^{\mu \ \alpha }_{\ \nu},
\ \Delta ^\mu _{\ \nu }(e^b-I), \ x^\alpha (e^b-I)$,
$\Delta ^\alpha _{\ \beta }\Delta ^\mu _{\ \nu },
\ x^{\alpha \beta }$.
Then we have the following:

{\bf Theorem 3.} ${\cal R}$ has the following properties:

(i)\ ${\cal R}$ is ad-invariant,

(ii)\ for any $a\in {\cal R}, \ S(a)^* \in {\cal R}$,

(iii)\ $ker\ \varepsilon /{\cal R}$ is spanned by the following
elements:

$$
x^\mu ;\ \ \ \Delta ^\mu _{\ \nu },\ \mu <\nu ;\ \ \ e^b-I .
$$

We see that our calculus is eleven-dimensional.
The basis of the space of the left-invariant 1-forms reads:

\beqn
\omega ^\mu  &=& \pi r^{-1}(I\otimes x^\mu ) =
e^{-b}\Lambda ^{\ \mu }_\nu dx^\nu , \nonumber \\
\omega ^\mu _{\ \nu } &=& \pi r^{-1}(I\otimes \Delta ^\mu _{\ \nu })
= \Lambda ^{\ \mu }_\alpha d\Lambda ^\alpha _{\ \nu }, \nonumber \\
\omega ^b &=& \pi r^{-1}(I\otimes (e^b-I)) = e^{-b}de^b. \label{fw}
\eeqn

The commutations rules between the invariant forms and generators of
${\cal W}_\kappa $ read:

\beqn
\lbrack e^b,\omega ^\mu   \rbrack  &=& 0,\ \ \ \ \
\lbrack  e^b,\omega ^\mu _{\ \nu } \rbrack  = 0,\ \ \ \ \
\lbrack e^b,\omega ^b  \rbrack  = 0,  \nonumber \\
\lbrack \Lambda ^\mu _{\ \nu },\omega ^b  \rbrack  &=& 0,\ \ \ \ \
\lbrack x^\alpha ,\omega ^b  \rbrack  = 0, \nonumber \\
\lbrack \Lambda ^\mu _{\ \nu },\omega ^\alpha   \rbrack  &=&
{i\over \kappa }(g_{\nu 0}\Lambda ^{\mu \alpha }
-\delta ^\alpha _{\ \nu }\Lambda ^\mu _{\ 0})\omega ^b   \nonumber \\
&\ &-{i\over \kappa }(g_{\nu 0}\Lambda ^\mu _{\ \tau }\omega ^{\tau \alpha }
+\Lambda ^\mu _{\ 0}\omega ^\alpha _{\ \nu }) , \nonumber \\
\lbrack x^\mu ,\omega ^\alpha   \rbrack  &=&
{i\over \kappa }e^b(\Lambda ^{\mu \alpha }\omega _0
-\delta ^\alpha _{\ 0}\Lambda ^\mu _{\ \nu }\omega ^\nu ), \nonumber \\
\lbrack \Lambda ^\alpha _{\ \beta },\omega ^\mu _{\ \nu }  \rbrack
 &=& 0, \nonumber \\
\lbrack  x^\alpha ,\omega ^\mu _{\ \nu } \rbrack  &=&
-{i\over \kappa }e^b(\delta ^\mu _{\ 0}\Lambda ^\alpha _{\ \beta }
\omega ^\beta _{\ \nu }
+g_{\nu 0}\Lambda ^\alpha _{\ \beta }\omega ^{\mu \beta }  \nonumber \\
&\ &-\Lambda ^\alpha _{\ \nu }\omega ^\mu _{\ 0}
-\Lambda ^{\alpha \mu }\omega _{0\nu }) . \label{kw}
\eeqn

The external power $\Gamma ^{\wedge 2}$
is described by the following relations:

\beqn
 &\ & \omega ^b \wedge \omega ^b = 0, \nonumber \\
 &\ & \omega ^\mu _{\ \nu }\wedge \omega ^\alpha _{\ \beta }
 +\omega ^\alpha _{\ \beta }\wedge \omega ^\mu _{\ \nu } = 0, \nonumber \\
 &\ & \omega ^b \wedge \omega ^\mu _{\ \nu }
 +\omega ^\mu _{\ \nu }\wedge \omega ^b = 0, \nonumber \\
 &\ & \omega ^b \wedge \omega ^\mu
 +\omega ^\mu \wedge \omega ^b = 0, \nonumber \\
 &\ & \omega ^\alpha \wedge \omega ^\mu
 +\omega ^\mu \wedge \omega ^\alpha
 +{i\over \kappa }(\delta ^\alpha _{\ 0}\omega ^\mu _{\ \beta }
 \wedge \omega ^\beta
 +\delta ^\mu _{\ 0} \omega ^\alpha _{\ \beta }
 \wedge \omega ^\beta  \nonumber \\
 &\ &  +\delta ^\mu _{\ 0}\omega ^b \wedge \omega ^\alpha
 +\delta ^\alpha _{\ 0}\omega ^b \wedge \omega ^\mu )
 +{2i\over \kappa }g^{\mu \alpha }\omega _0 \wedge \omega ^b =0, \nonumber \\
 &\ & \omega ^\alpha \wedge \omega ^{\mu \nu }
 +\omega ^{\mu \nu }\wedge \omega ^\alpha
 +{i\over \kappa }(\delta ^\mu _{\ 0}\omega ^{\tau \nu }
 \wedge \omega ^{\ \alpha }_\tau
 +\delta ^\nu _{\ 0}\omega ^\mu _{\ \tau }\wedge
 \omega ^{\tau \alpha }   \nonumber \\
&\ & +\omega ^{\ \nu }_0 \wedge\omega ^{\alpha \mu }
 +\omega ^\mu _{\ 0}\wedge \omega ^{\alpha \nu })
  +{i\over \kappa }(\delta ^\mu _{\ 0}\omega ^b \wedge
 \omega ^{\alpha \nu }
 +\delta ^\nu _{\ 0}\omega ^b \wedge \omega ^{\mu \alpha }  \nonumber \\
&\ & +g^{\alpha \mu }\omega ^{\ \nu }_0 \wedge \omega ^b
 +g^{\alpha \nu }\omega ^\mu _{\ 0}\wedge \omega ^b) = 0.  \label{zw}
\eeqn

The basis of $\Gamma ^{\wedge 2}$ consists of the following elements:

\beqn
 &\ & \omega ^{\alpha \beta }\wedge \omega ^{\mu \nu };
 \ \ \ \ {\rm for}\ \alpha <\beta ,\ \mu < \nu,\ (\alpha \beta )
 \ne (\mu \nu ),\ \alpha < \mu ; \nonumber \\
 &\ & \omega ^{\mu \nu }\wedge \omega ^b,
 \ \omega ^\mu \wedge \omega ^\nu ,
 \ \omega ^{\mu \nu }\wedge \omega ^\alpha ;
 \ \ \ \ {\rm for}\ \mu <\nu ; \nonumber \\
 &\ & \omega ^\mu \wedge \omega ^b . \nonumber
\eeqn

The Cartan-Maurer equations have the following form:
\beqn
d\omega ^\mu _{\ \nu } &=& \omega ^{\ \mu }_\tau
\wedge \omega ^\tau _{\ \nu }, \nonumber \\
d\omega ^\mu  &=& \omega ^\mu \wedge \omega ^b
+ \omega ^{\ \mu }_\tau \wedge \omega ^\tau , \nonumber \\
d\omega ^b &=& 0. \label{mw}
\eeqn

Note that if we put $b=0$ in eq.(\ref{fw}), (\ref{kw}),
we obtain the same sets of relations like
in eq.(\ref{fp}), (\ref{kp}),
 but it is not true for
the case of eq.(\ref{zw}), (\ref{mw})
end eq.(\ref{zp}), (\ref{mp}).

In order to obtain the quantum Lie algebra we introduce
the left-invariant field defined by the formula:

$$
da = {1\over 2}(\chi _{\mu \nu }*a)\omega ^{\mu \nu }
+(\chi _\mu *a)\omega ^\mu
+(\lambda *a)\omega ^b .
$$

The resulting quantum Lie algebra reads:
\beqn
\lbrack \chi _{\mu \nu },\chi _\alpha  \rbrack &=&
(1+{i\over \kappa }\chi _0)
(\chi _\mu g_{\alpha \nu }-\chi _\nu g_{\alpha \mu }) ,\nonumber  \\
\lbrack \chi _\alpha ,\chi _\mu   \rbrack &=& 0   \nonumber \\
\lbrack \chi _\mu ,\lambda   \rbrack &=&  -\chi _\mu
-{i\over \kappa }\chi _0 \chi _\mu
+{i\over \kappa }g_{0\mu }g^{\sigma \alpha }\chi _\sigma
\chi _\alpha , \nonumber \\
\lbrack \chi _{\mu \nu },\lambda   \rbrack &=& -{i\over \kappa }
(\chi _\mu \chi _{0\nu }+\chi _\nu \chi _{\mu 0})  \nonumber \\
&\ &+{i\over \kappa }g^{\alpha \sigma }(\chi _\alpha
\chi _{\sigma \nu }g_{0\mu }
+\chi _\alpha \chi _{\mu \sigma }g_{0\nu }), \nonumber  \\
\lbrack \chi _{\alpha \beta },\chi _{\mu \nu } \rbrack &=&
(\chi _{\beta \mu }g_{\nu \alpha }
+\chi _{\alpha \nu }g_{\mu \beta }
-\chi _{\beta \nu }g_{\mu \alpha }
-\chi _{\alpha \mu }g_{\nu \beta }) \nonumber \\
 +{i\over \kappa }(\chi _\alpha (
\chi _{\mu 0}g_{\nu \beta }&-&\chi _{\nu 0}g_{\mu \beta })
+\chi _\mu (\chi _{\beta 0}g_{\alpha \nu }
-\chi _{\alpha 0}g_{\nu \beta }) \nonumber \\
 +\chi _\nu (\chi _{\alpha 0}g_{\mu \beta }
&-&\chi _{\beta 0}g_{\alpha \mu })
+\chi _\beta (\chi _{\nu 0}g_{\alpha \mu }
-\chi _{\mu 0}g_{\alpha \nu })) \nonumber \\
 +{i\over \kappa }(\chi _\beta (\chi _{\mu \alpha }g_{0\nu }
&-&\chi _{\nu \alpha }g_{0\mu })
+\chi _\nu (\chi _{\beta \mu }g_{0\alpha }
-\chi _{\alpha \mu }g_{0\beta }) \nonumber \\
 +\chi _\mu (\chi _{\alpha \nu }g_{0\beta }
&-&\chi _{\beta \nu }g_{0\alpha })
+\chi _\alpha (\chi _{\nu \beta }g_{0\mu }
-\chi _{\mu \beta }g_{0\nu })) .\nonumber
\eeqn

The Weyl algebra reads~\cite{defwayla}:

{\noindent The commutation rules:}
\beqn
\lbrack  M^{ij},P_0 \rbrack  &=& 0 \nonumber \\
\lbrack M^{ij},P_k  \rbrack  &=& i\kappa (\delta ^j_{\ k}
g^{0i}-\delta ^i_{\ k}g^{0j})(1-e^{-{P_0\over \kappa }})
+ i(\delta ^j_{\ k}g^{is}-\delta ^i_{\ k}g^{js})P_s\nonumber \\
\lbrack M^{i0},P_0  \rbrack  &=& i\kappa g^{i0}
(1-e^{-{p_0\over \kappa }})+ig^{ik}P_k \nonumber \\
\lbrack M^{i0},P_k  \rbrack  &=& -i{\kappa \over 2}g^{00}
\delta ^i_{\ k}(1-e^{-2{P_0\over \kappa }})
-i\delta ^i_{\ k}g^{0s}P_s e^{-{P_0\over \kappa }}+ \nonumber \\
 &\ & +ig^{0i}P_k (e^{-{P_0\over \kappa }}-1)+
{i\over 2\kappa }\delta ^i_{\ k}g^{rs}P_r P_s -
{i\over \kappa }g^{is} P_s P_k \nonumber \\
\lbrack P_\mu ,P_\nu   \rbrack  &=& 0 ,\nonumber \\
\lbrack M^{\mu \nu },M^{\lambda \sigma }  \rbrack
&=& i(g^{\mu \sigma }M^{\nu \lambda }-g^{\nu \sigma }M^{\mu \lambda }
+g^{\nu \lambda }M^{\mu \sigma }-g^{\mu \lambda }M^{\nu \sigma }), \nonumber \\
\lbrack M_{\mu \nu },D  \rbrack  &=& 0, \nonumber \\
\lbrack D,P_0  \rbrack  &=& i\kappa (1-e^{-{P_0\over \kappa }}) \nonumber \\
\lbrack D,P_i  \rbrack  &=& iP_i e^{-{P_0\over \kappa }}
+i{\kappa \over 2}g^{00}g_{i0}(1-e^{-{P_0\over \kappa }})^2 +  \nonumber \\
 &\ & +ig_{0i}g^{0s}P_s (1-e^{-{P_0\over \kappa }})
+{i\over 2\kappa }g_{0i}g^{rs}P_r P_s \nonumber
\eeqn
The coproducts, counit and antipode:
\beqn
\Delta D &=& D\otimes I + I\otimes D -
g_{0i}M^{i0}\otimes (1-e^{-{P_0\over \kappa }})-
{1\over \kappa }g_{0i}M^{ik}\otimes P_k \nonumber \\
\Delta P_0 &=& I\otimes P_0 + P_0 \otimes I \nonumber \\
\Delta P_k &=& P_k \otimes e^{-{P_0\over \kappa }}+
I\otimes P_k \nonumber \\
\Delta M^{ij} &=& M^{ij}\otimes I + I\otimes M^{ij} \nonumber \\
\Delta M^{i0} &=& I\otimes M^{i0} + M^{i0}\otimes
e^{-{P_0\over \kappa }}- {1\over \kappa }M^{ij}\otimes P_j , \nonumber \\
\varepsilon (M^{\mu \nu }) &=& 0; \ \ \
\varepsilon (P_\nu ) = 0;\ \ \
\varepsilon (D) =0, \nonumber \\
S(P_0) &=& -P_0 ,\nonumber \\
S(P_i) &=& -e^{P_0\over \kappa }P_i ,\nonumber \\
S(M^{ij}) &=& -M^{ij} ,\nonumber \\
S(M^{i0}) &=& -e^{P_0\over \kappa }(M^{i0}+
{1\over \kappa }M^{ij}P_j) ,\nonumber \\
S(D) &=& -D+g_{i0}M^{i0}-g_{i0}e^{P_0\over \kappa }M^{i0}
-{1\over \kappa }g_{i0}e^{P_0\over \kappa }M^{ik}P_k .\nonumber
\eeqn

Using, on the one hand, the properties of the left-invariant fields
described by Woronowicz~\cite{rach} and on the other hand,
the duality relations ${\cal W}_\kappa \iff \tilde {\cal W}_\kappa $
established in~\cite{w11}, one can prove that the following substitutions:

\beqn
\chi _0  &=& -i\kappa (e^{P_0\over \kappa }-1) ,\nonumber \\
\chi _i  &=& -i(e^{P_0\over \kappa }P_i
-g_{i0}{M^2 \over 2\kappa }) ,\nonumber \\
\chi _{ij} &=& -i (1+{i\over \kappa }\chi _0 )M_{ij}
-{1\over \kappa }(\chi _0 M_{ij} +\chi _i M_{j0}- \chi _j M_{i0}) ,\nonumber \\
\chi _{i0} &=& -i (1+{i\over \kappa }\chi _0 )M_{i0} ,\nonumber \\
\lambda  &=& -iD -{1\over \kappa }\chi ^\alpha M_{\alpha 0} ,\nonumber
\eeqn
where
\beqn
M^2 &=& g^{00}(2\kappa \sinh({P_0\over 2\kappa }))^2
+4\kappa g^{0l} P_l e^{P_0\over 2\kappa }
\sinh({P_0\over 2\kappa })  \nonumber \\
&\ &+g^{rs}P_r e^{P_0\over 2\kappa }P_s e^{P_0\over 2\kappa }, \nonumber
\eeqn

reproduce the algebra and coalgebra structure of our quantum Lie algebra.

\section{Acknowledgment}
\zero

The numerous discussions with  dr Pawe\l{}  Ma\`slanka are kindly acknowledged.


\begin{thebibliography}{99}
\newcommand{\by}[1]{#1 }
\newcommand{\jo}[1]{ #1 }
\newcommand{\ty}[1]{{\it #1} }
\newcommand{\vol}[1]{{\bf #1} }
\newcommand{\ef}{.}
\newcommand{\yr}[1]{(#1), }

\bibitem{defwayla}\by{P.Kosi\'nski, P.Ma\'slanka}
\ty{The $\kappa -$Weyl group and its algebra}
\jo{in "From Field Theory to Quantum Groups"}vol. on 60
anniversary of J.Lukierski, World Scientific Singapur, 1996
or \vol{Q-ALG}9512018\ef

%
\bibitem{diag}\by{P.Kosi\'nski, P.Ma\'slanka, Jan Sobczyk}
\ty{The bicovariant diffrerential calculus on the $\kappa -$Poincar\`e
group and on the $\kappa -$Minkowski space}
\jo{Czechoslovak J. Phys.}\vol{46}\yr{1996}201-208 or \vol{Q-ALG}9508021\ef

\bibitem{tree}\by{P.Kosi\`nski, M.Majewski, P.Ma\`slanka}
\ty{The bicovariant differential calculus on the three-dimensional
Kappa-Poincar'e group}\vol{Q-ALG}9602004,
in print in Acta Phys. Pol. \vol{B}\ef

\bibitem{min}\by{C.Gonera, P.Kosi\`nski, P.Ma\`slanka}
\ty{Differential calculi on quantum Minkowski space}
\vol{Q-ALG}9602007, in print in J. Math. Phys\ef


\bibitem{rach}\by{S.L.Woronowicz}
\ty{Differential Calculus on Compact Ma\-trix Pseudo\-groups
(Quantums Groups)}\jo{Comm. Math. Phys.}\vol{122}\yr{1989}125\ef

\bibitem{luk}\by{P.Kosi\'nski, J.Lukierski, P.Ma\'slanka, A.Sitarz}
-paper in preparation\ef

\bibitem{w1}\by{W.B.Schmidke, J.Weiss, B.Zumino}\jo{Zeitschr. f. Physik}
\vol{52}\yr{1991}472\ef

\by{U.Carow-Watamura, M.Schliecker, M.Scholl, S.Watamura}
\jo{Int. J. Mod. Phys.}\vol{A\,6}\yr{1991}3081\ef

\by{S.L. Woronowicz}\jo{Comm. Math. Phys.}\vol{136}\yr{1991}399\ef

\by{O. Ogievetsky, W.B. Schmidke, J. Weiss, B. Zumino}
\jo{Comm. Math.  Phys.}\vol{150}\yr{1992}495\ef

\by{M. Chaichian, A.P. Demichev}\jo{Proceedings of the
Workshop: "Generalized symmetries in Physics", Clausthal}\yr{1993}\ef

\by{V. Dobrev}\jo{J.  Phys.}\vol{A\,26}\yr{1993}1317\ef

\by{L. Castellani}\jo{in "Quantum Groups" Proceedings of XXX
Karpacz Winter School of Theoretical Physics, Karpacz 1994, PWN 1995, p.
13}\ef

\bibitem{w2}\by{J. Lukierski, A. Nowicki, H. Ruegg, V. Tolstoy}\jo{Phys.
Lett.}\vol{B\,264}\yr{1991}331\ef

\by{J. Lukierski, A. Nowicki, H. Ruegg}\jo{Phys.
Lett.}\vol{B\,293}\yr{1993}344\ef

\by{S. Giller,  P. Kosi\'nski, J. Kunz,  M. Majewski, P.
Ma\'slanka}\jo{Phys. Lett.}\vol{B\,286}\yr{1992}57\ef

\bibitem{w3}\by{S. Zakrzewski }\jo{J. Phys.}\vol{A\,27}\yr{1994}2075\ef

\bibitem{w4}\by{ J. Lukierski, H. Ruegg }\jo{Phys. Lett.}
\vol{B\,329}\yr{1994}189\ef

\by{P. Ma\'slanka}\jo{J. Phys.}\vol{A\,26}\yr{1993}L1251\ef

\bibitem{w5}\by{S. Giller,  P. Kosi\'nski, J. Kunz,  M. Majewski, P.
Ma\'slanka}\jo{Phys. Lett.}\vol{B\,286}\yr{1992}57\ef

\by{S. Giller, C. Gonera,  P. Kosi\'nski, J. Kunz,  P.
Ma\'slanka}\jo{Mod. Phys. Lett.}\vol{A\,8}\yr{1993}3785\ef

\by{ J. Lukierski, H. Ruegg, W. R\"uhl }\jo{Phys. Lett.}\vol{B\,313}
\yr{1993}357\ef

\bibitem{w6}\by{P. Ma\'slanka}\jo{J. Math. Phys.}\vol{34}\yr{1993}
6025\ef

\bibitem{w7}\by{---------}\jo{J.  Phys. }\vol{A\,27}\yr{1994}7099\ef

\bibitem{w8}\by{S. Majid}\jo{Pacific J. Math.}\vol{141}\yr{1990}311\ef

\by{---------}\jo{Int. J. Mod. Phys.}\vol{A\,5}\yr{1990}1\ef

\by{---------}\jo{J. Algebra}\vol{130}\yr{1990}17\ef

\bibitem{w9}\by{---------, H. Ruegg}\jo{Phys. Lett.}\vol{B\,334}
\yr{1994}348\ef

\by{Ph. Zaugg}\jo{preprint MIT--CTP--2353}\yr{1994}\ef

\bibitem{w10}\by{A. Ballesteros, E. Celeghini, R. Giachetti, E.
Sorace, M. Tarlini}\jo{J. Phys.}\vol{A\,26}\yr{1993}7495\ef

\by{P. Ma\'slanka}\jo{J. Math. Phys. }\vol{35}\yr{1994}1976\ef

\by{Ph. Zaugg}\jo{preprint MIT--CTP--2353}\yr{1994}\ef

\bibitem{w11}\by{P. Kosi\'nski,  P. Ma\'slanka}\ty{preprint IMU\L \/
3\slash 94}or \vol{Q-ALG}9411033\ef


\end{thebibliography}
\end{document}